\documentclass[11pt]{amsart}

\theoremstyle{plain}
\newtheorem{theorem}                 {Theorem}      [section]

\newtheorem{corollary}    [theorem]  {Corollary}
\newtheorem{lemma}        [theorem]  {Lemma}

\theoremstyle{definition}

\DeclareMathOperator{\Rim}{Rim}

\title{Conformally flat tangent bundles with general
natural lifted metrics}

\author{S.~L.~Dru\c t\u a*}

\thanks{*  Partially supported by the Grant ET 5871; 2006,2007,
CEEX, Ministerul Educa\c tiei \c si Cercet\u arii, Rom\^ania}

\begin{document}
\maketitle

\maketitle

\begin{minipage}{2.5in}
\begin{flushleft}
Contemporary Geometry and Topology and Related Topics Cluj-Napoca,
August 19-25, \\
 2007, pp. 153–-166
 \end{flushleft}
\end{minipage}

\begin{abstract}

We study the conditions under which the tangent bundle $(TM,G)$ of
an $n$-dimensional Riemannian manifold $(M,g)$ is conformally flat,
where $G$ is a general natural lifted metric of $g$. We prove that
the base manifold must have constant sectional curvature and we find
some expressions for the natural lifted metric $G$, such that the
tangent bundle $(TM,G)$ become conformally flat.

{\it Mathematics Subject Classification 2000:} Primary 53C55, 53C15, 53C05\\
\rightskip=1.2cm \leftskip=1.2cm {\it Key words and phrases}:
tangent bundle, Riemannian metric, general natural lift, conformal
curvature.
\end{abstract}

\section{Introduction}

The possibility to consider vertical, complete and horizontal
lifts on the tangent bundle $TM$ of a smooth $n$-dimensional
Riemannian manifold $(M,g)$, leads to some interesting geometric
structures, studied in the last years (see \cite{Abbassi1},
\cite{Abbassi2}, \cite{BejanOpr}, \cite{Munteanu1},
\cite{Munteanu2}), and to interesting relations with some problems
in Lagrangian and Hamiltonian mechanics. One uses several
Riemannian and pseudo-Riemannian metrics, induced by the
Riemannian metric $g$ on $M$. Among them, we may quote the Sasaki
metric, the Cheeger-Gromoll metric and the complete lift of the
metric $g$. On the other hand, the natural lifts of $g$ to $TM$,
introduced in the papers \cite{KowalskiSek} and \cite{Krupka},
induce some new Riemannian and pseudo-Riemannian geometric
structures with many nice geometric properties (\cite{Kolar},
\cite{KowalskiSek}).

Professor Oproiu has studied some properties of a natural lift
$G$, of diagonal type, of the Riemannian metric $g$ and a natural
almost complex structure $J$ of diagonal type  on $TM$ (see
\cite{Oproiu1}, \cite{Oproiu2}, \cite{Oproiu3}, and see also
\cite{OprPap1}, \cite{OprPap2}). In the paper \cite{Oproiu4}, the
same author has presented a general expression of the natural
almost complex structures on $TM$. In the definition of the
natural almost complex structure $J$ of general type there are
involved eight parameters (smooth functions of the density energy
on $TM$). However, from the condition for $J$ to define an almost
complex structure, four of the above parameters can be expressed
as (rational) functions of the other four parameters. A Riemannian
metric $G$ which is a natural lift of general type of the metric
$g$ depends on other six parameters.

In \cite{Druta}, the present author got the conditions and the
unique form of the matrix associated to the general natural lifted
metric $G$, such that the the tangent bundle $TM$, with respect to
the metric $G$ has constant sectional curvature. In
\cite{OprDruta} we have found the conditions under which the
K\"ahlerian manifold $(TM,G,J)$ has constant holomorphic sectional
curvature.

In the present paper we study the conformal curvature of the
tangent bundle of a Riemannian manifold $(M,g)$. Namely, we are
interested in finding the conditions under which the Riemannian
manifold $(TM,G)$, where $G$ is the general natural lifted metric
of $g$, is conformally flat.

\section{Preliminary results}

Consider a smooth $n$-dimensional Riemannian manifold $(M,g)$ and
denote its tangent bundle by $\tau :TM\longrightarrow M$. Recall
that $TM$ has a structure of a $2n$-dimensional smooth manifold,
induced from the smooth manifold structure of $M$. This structure
is obtained by using local charts on $TM$ induced  from usual
local charts on $M$. If $(U,\varphi )= (U,x^1,\dots ,x^n)$ is a
local chart on $M$, then the corresponding induced local chart on
$TM$ is $(\tau ^{-1}(U),\Phi )=(\tau ^{-1}(U),x^1,\dots , x^n,$
$y^1,\dots ,y^n)$, where the local coordinates $x^i,y^j,\
i,j=1,\dots ,n$, are defined as follows. The first $n$ local
coordinates of a tangent vector $y\in \tau ^{-1}(U)$ are the local
coordinates in the local chart $(U,\varphi)$ of its base point,
i.e. $x^i=x^i\circ \tau$, by an abuse of notation. The last $n$
local coordinates $y^j,\ j=1,\dots ,n$, of $y\in \tau ^{-1}(U)$
are the vector space coordinates of $y$ with respect to the
natural basis in $T_{\tau(y)}M$ defined by the local chart
$(U,\varphi )$. Due to this special structure of differentiable
manifold for $TM$, it is possible to introduce the concept of
$M$-tensor field on it (see \cite{Mok}). The $M$-tensor fields are
defined by their components with respect to the induced local
charts on $TM$ (hence they are defined locally), but they can be
interpreted as some (partial) usual tensor fields on $TM$.
However, the essential quality of an $M$-tensor field on $TM$ is
that the local coordinate change rule of its components with
respect to the change of induced local charts is the same as the
local coordinate change rule of the components of an usual tensor
field on $M$ with respect to the change of local charts on $M$.
More precisely, an $M$-tensor field of type $(p,q)$ on $TM$ is
defined by sets of $n^{p+q}$ components (functions depending on
$x^i$ and $y^i$), with $p$ upper indices and $q$ lower indices,
assigned to induced local charts $(\tau ^{-1}(U),\Phi )$ on $TM$,
such that the local coordinate change rule of these components
(with respect to induced local charts on $TM$) is that of the
local coordinate components  of a tensor field of type $(p,q)$ on
the base manifold $M$ (with respect to usual local charts on $M$),
when a change of local charts on $M$ (and hence on $TM$) is
performed (see \cite{Mok} for further details); e.g., the
components $y^i,\ i=1,\dots ,n$, corresponding to the last $n$
local coordinates of a tangent vector $y$, assigned to the induced
local chart $(\tau ^{-1}(U), \Phi )$ define an $M$-tensor field of
type $(1,0)$ on $TM$. A usual tensor field of type $(p,q)$ on $M$
may be thought of as an $M$-tensor field of type $(p,q)$ on $TM$.
If the considered tensor field on $M$ is covariant only, the
corresponding $M$-tensor field on $TM$ may be identified with the
induced (pullback by $\tau $) tensor field on $TM$.  Some useful
$M$-tensor fields on $TM$ may be obtained as follows. Let
$u:[0,\infty )\longrightarrow {\mathbb{R}}$ be a smooth function
and let $\|y\|^2=g_{\tau (y)}(y,y)$ be the square of the norm of
the tangent vector $y\in \tau ^{-1}(U)$. If $\delta ^i_j$ are the
Kronecker symbols (in fact, they are the local coordinate
components of the identity tensor field $I$  on $M$), then the
components $u(\|y\|^2)\delta ^i_j$ define an $M$-tensor field of
type $(1,1)$ on $TM$. Similarly, if $g_{ij}(x)$ are the local
coordinate components of the metric tensor field $g$ on $M$ in the
local chart $(U,\varphi )$, then the components $u(\|y\|^2)
g_{ij}$ define a symmetric $M$-tensor field of type $(0,2)$ on
$TM$. The components $g_{0i}=y^kg_{ki}$ define an $M$-tensor field
of type $(0,1)$ on $TM$.

Denote by $\dot \nabla$ the Levi Civita connection of the
Riemannian metric $g$ on $M$. Then we have the direct sum
decomposition
\begin{equation}
TTM=VTM\oplus HTM
\end{equation}
of the tangent bundle to $TM$ into the vertical distribution
$VTM={\rm Ker}\ \tau_*$ and the horizontal distribution $HTM$
defined by $\dot \nabla $. The set of vector fields
$(\frac{\partial}{\partial y^1}, \dots , \frac{\partial}{\partial
y^n})$ on $\tau ^{-1}(U)$ defines a local frame field for $VTM$
and for $HTM$ we have the local frame field $(\frac{\delta}{\delta
x^1},\dots ,\frac{\delta}{\delta x^n})$, where
$$
\frac{\delta}{\delta x^i}=\frac{\partial}{\partial
x^i}-\Gamma^h_{0i} \frac{\partial}{\partial y^h},\ \ \ \Gamma
^h_{0i}=y^k\Gamma ^h_{ki},
 $$
and $\Gamma ^h_{ki}(x)$ are the Christoffel symbols of $g$.

The set $(\frac{\partial}{\partial y^1},\dots
,\frac{\partial}{\partial y^n}, \frac{\delta}{\delta x^1},\dots
,\frac{\delta}{\delta x^n})$ defines a local frame on $TM$,
adapted to the direct sum decomposition (1). Remark that
$$
\frac{\partial}{\partial y^i}=(\frac{\partial}{\partial x^i})^V,\
\ \frac{\delta}{\delta x^i}=(\frac{\partial}{\partial x^i})^H,
$$
where $X^V$ and $X^H$ denote the vertical and horizontal lift of
the vector field $X$ on $M$ respectively. We can use the vertical
and horizontal lifts  in order to obtain invariant expressions for
some results in this paper. However, we should prefer to work in
local coordinates  since the formulas are obtained easier and, in
a certain sense, they are more natural.

We can easily obtain the following

\begin{lemma}\label{lema1}
If $n>1$ and $u,v$ are smooth functions on $TM$ such that
$$
u g_{ij}+v g_{0i}g_{0j}=0,
$$
on the domain of any induced local chart on $TM$, then $u=0,\
v=0$.
\end{lemma}

\bf Remark. \rm In a similar way  we obtain from the condition
$$
u\delta ^i_j+vg_{0j}y^i=0
$$
the relations $u=v=0$.

Consider the energy density of the tangent vector $y$ with respect
to the Riemannian metric $g$
\begin{equation}
t=\frac{1}{2}\|y\|^2=\frac{1}{2}g_{\tau(y)}(y,y)=\frac{1}{2}g_{ik}(x)y^iy^k,
\ \ \ y\in \tau^{-1}(U).
\end{equation}
Obviously, we have $t\in [0,\infty)$ for all $y\in TM$.

\section{The conformal curvature of the tangent bundle with general natural lifted metric}

Let $G$ be the general natural lifted metric on $TM$, defined by
\begin{equation}\label{defG}
\begin{array}{l}
G(\frac{\delta}{\delta x^i}, \frac{\delta}{\delta x^j})=c_1g_{ij}+
d_1g_{0i}g_{0j}=G^{(1)}_{ij}
\\   \mbox{ } \\
G(\frac{\partial}{\partial y^i}, \frac{\partial}{\partial y^j})=
c_2g_{ij}+d_2g_{0i}g_{0j}=G^{(2)}_{ij}
\\   \mbox{ } \\
G(\frac{\partial}{\partial y^i},\frac{\delta}{\delta x^j})=
G(\frac{\delta}{\delta x^i},\frac{\partial}{\partial y^j})=
c_3g_{ij}+d_3g_{0i}g_{0j}=G^{(3)}_{ij},
\end{array}
\end{equation}
where $c_1,c_2,c_3,d_1,d_2,d_3$ are six smooth functions of the
density energy on $TM$.

The Levi-Civita connection $\nabla$ of the Riemannian manifold
$(TM,G)$ is obtained from the formula
$$
2G(\nabla_XY,Z)=X(G(X,Z))+Y(G(X,Z))-Z(G(X,Y))
$$
$$
+G([X,Y],Z)-G([X,Z],Y)-G([Y,Z],X); ~~\forall X,Y,Z\in \chi(M)
$$
and is characterized by the conditions
$$
\nabla G=0,\ T=0,
$$
where $T$ is the torsion tensor of $\nabla.$

In the case of the tangent bundle $TM$ we can obtain the explicit
expression of $\nabla$. The symmetric $2n\times 2n$ matrix
$$
\left(
\begin{array}{ll}
G^{(1)}_{ij} &  G^{(3)}_{ij} \\
G^{(3)}_{ij} &  G^{(2)}_{ij}
\end{array}
\right)
$$
associated to the metric $G$ in the base $(\frac{\delta}{\delta
x^1},\dots ,\frac{\delta}{\delta x^n},\frac{\partial}{\partial
y^1},\dots ,\frac{\partial}{\partial y^n})$ has the inverse
$$
\left(
\begin{array}{ll}
H_{(1)}^{ij} &  H_{(3)}^{ij} \\
H_{(3)}^{ij} &  H_{(2)}^{ij}
\end{array}
\right)
$$
where the entries are the blocks
 \begin{eqnarray*}
H_{(1)}^{kl}=p_1g^{kl}+q_1y^ky^l
\end{eqnarray*}
\begin{eqnarray}\label{matrinv}
H_{(2)}^{kl}=p_2g^{kl}+q_2y^ky^l
\end{eqnarray}
\begin{eqnarray*}
H_{(3)}^{kl}=p_3g^{kl}+q_3y^ky^l.
\end{eqnarray*}
Here $g^{kl}$ are the components of the inverse of the matrix
$(g_{ij})$ and $p_1,q_1,p_2,q_2,p_3$, $q_3:[0,\infty)\rightarrow
\mathbb{R},$ some real smooth functions. Their expressions are
obtained  by solving the system:
$$
\begin{cases}
G^{(1)}_{ih}H_{(1)}^{hk}+G^{(3)}_{ih}H_{(3)}^{hk}=\delta_i^k\\
G^{(1)}_{ih}H_{(3)}^{hk}+G^{(3)}_{ih}H_{(2)}^{hk}=0\\
G^{(3)}_{ih}H_{(1)}^{hk}+G^{(2)}_{ih}H_{(3)}^{hk}=0\\
G^{(3)}_{ih}H_{(3)}^{hk}+G^{(2)}_{ih}H_{(2)}^{hk}=\delta_i^k,
\end{cases}
$$
in which we substitute the relations (\ref{defG}) and
(\ref{matrinv}). By using Lemma \ref{lema1},  we get $p_1,p_2,p_3$
as functions of $c_1,c_2,c_3$

\begin{eqnarray}\label{inversa1}
p_1=\frac{c_2}{c_1c_2-c_3^2},\ \ p_2=\frac{c_1}{c_1c_2-c_3^2},\ \
p_3=-\frac{c_3}{c_1c_2-c_3^2}
\end{eqnarray}
and $q_1,q_2,q_3$ as functions of $c_1,c_2,c_3,$ $d_1,d_2,d_3,$
$p_1,p_2,p_3$
\begin{eqnarray*}
q_1=-\frac{c_2d_1p_1 - c_3d_3p_1 - c_3d_2p_3 + c_2d_3p_3 +
2d_1d_2p_1t - 2d_3^2p_1t}{c_1c_2 - c_3^2 + 2c_2d_1t + 2c_1d_2t -
4c_3d_3t + 4d_1d_2t^2 - 4d_3^2t^2},
\end{eqnarray*}
\begin{eqnarray}\label{inversa2}
q_2=-\frac{d_2p_2 + d_3p_3}{c_2 + 2d_2t} +
\end{eqnarray}
\begin{eqnarray*}
+ \frac{(c_3 + 2d_3t)[(d_3p_1 + d_2p_3)(c_1 + 2d_1t) - (d_1p_1 +
d_3p_3)(c_3 + 2d_3t)]}{(c_2 + 2d_2t)[(c_1 + 2d_1t)(c_2 + 2d_2t) -
(c_3 + 2d_3t)^2]}
\end{eqnarray*}
\begin{eqnarray*}
q_3=-\frac{(d_3p_1 + d_2p_3)(c_1 + 2d_1t) - (d_1p1 + d_3p_3)(c_3 +
2d_3t)}{(c_1 + 2d_1t)(c_2 + 2d_2t) - (c_3 + 2d_3t)^2}.
\end{eqnarray*}
In the paper \cite{OprDruta} we obtained the expression of the
Levi Civita connection of the Riemannian metric $G$ on $TM$.
\begin{theorem}
The Levi-Civita connection $\nabla$ of\ $G$ has the following
expression in the local adapted frame $(\frac{\partial}{\partial
y^1}, \dots , \frac{\partial}{\partial y^n},\frac{\delta}{\delta
x^1},\dots ,\frac{\delta}{\delta x^n})$
$$
\begin{cases}
$$
\nabla_{\frac{\partial}{\partial y^i}} \frac{\partial}{\partial
y^j}=Q^{h}_{ij}\frac{\partial}{\partial
y^h}+\widetilde{Q}^{h}_{ij}\frac{\delta}{\delta x^h},~
\nabla_{\frac{\delta}{\delta x^i}} \frac{\partial}{\partial
y^j}=(\Gamma^h_{ij}+\widetilde{P}^{h}_{ji})\frac{\partial}{\partial
y^h}+P^{h}_{ji}\frac{\delta}{\delta x^h}\\

\nabla_{\frac{\partial}{\partial y^i}} \frac{\delta}{\delta
x^j}=P^{h}_{ij}\frac{\delta}{\delta
x^h}+\widetilde{P}^{h}_{ij}\frac{\partial}{\partial y^h},~
\nabla_{\frac{\delta}{\delta x^i}} \frac{\delta}{\delta
x^j}=(\Gamma^h_{ij}+\widetilde{S}^{h}_{ij})\frac{\delta}{\delta
x^h}+S^{h}_{ij}\frac{\partial}{\partial y^h},
$$
\end{cases}
$$
where $\Gamma^h_{ij}$ are the Christoffel symbols of the
connection $\dot\nabla$ and the $M$-tensor fields appearing as
coefficients in the above expressions are given as
$$
\left\{
\begin{array}{l}
Q_{ij}^h=\frac{1}{2}(\partial_iG_{jk}^{(2)}+\partial_jG_{ik}^{(2)}-
\partial_kG_{ij}^{(2)})H_{(2)}^{kh}+\frac{1}{2}(\partial_iG_{jk}^{(3)}+
\partial_jG_{ik}^{(3)})H_{(3)}^{kh}\\
\widetilde{Q}_{ij}^h=\frac{1}{2}(\partial_iG_{jk}^{(2)}+\partial_jG_{ik}^{(2)}-
\partial_kG_{ij}^{(2)})H_{(3)}^{kh}+\frac{1}{2}(\partial_iG_{jk}^{(3)}+
\partial_jG_{ik}^{(3)})H_{(1)}^{kh}\\\\
P^h_{ij}=\frac{1}{2}(\partial_iG_{jk}^{(3)}-
\partial_kG_{ij}^{(3)})H_{(3)}^{kh}+\frac{1}{2}(\partial_iG_{jk}^{(1)}+
R^l_{0jk}G^{(2)}_{li})H_{(1)}^{kh}\\
\widetilde{P}^h_{ij}=\frac{1}{2}(\partial_iG_{jk}^{(3)}-
\partial_kG_{ij}^{(3)})H_{(2)}^{kh}+\frac{1}{2}(\partial_iG_{jk}^{(1)}+
R^l_{0jk}G^{(2)}_{li})H_{(3)}^{kh}\\\\
S^h_{ij}=-\frac{1}{2}(\partial_kG_{ij}^{(2)}+R^l_{0ij}G^{(2)}_{lk})H_{(2)}^{kh}+c_3R_{i0jk}H_{(3)}^{kh}\\
\widetilde{S}^h_{ij}=-\frac{1}{2}(\partial_kG_{ij}^{(1)}+R^l_{0ij}G^{(2)}_{lk})H_{(3)}^{kh}+c_3R_{i0jk}H_{(1)}^{kh},
\end{array}
\right.
$$
where $R^h_{kij}$ are the components of the curvature tensor field
of the Levi Civita connection $\dot \nabla$ of the base manifold
$(M,g)$.
\end{theorem}

Taking into account the expressions (\ref{defG}), (\ref{matrinv})
 and by using the formulas (\ref{inversa1}), (\ref{inversa2}) we
can obtain the detailed expressions of
$P^h_{ij},Q^h_{ij},S^h_{ij}, \widetilde P^h_{ij}, \widetilde
Q^h_{ij}, \widetilde S^h_{ij}.$

The curvature tensor field $K$ of the connection $\nabla$ is
defined by the well known formula
$$
K(X,Y)Z=\nabla_X\nabla_YZ-\nabla_Y\nabla_XZ-\nabla_{[X,Y]}Z,\ \ \
X,Y,Z\in \Gamma (TM).
$$

By using the local adapted frame $(\frac{\delta}{\delta
x^i},\frac{\partial}{\partial y^j}),\ i,j=1,\dots ,n$ we obtained
in \cite{OprDruta}, after a standard straightforward computation
$$
K\big(\frac{\delta}{\delta x^i},\frac{\delta}{\delta x^j}
\big)\frac{\delta}{\delta x^k}=XXXX^h_{kij}\frac{\delta}{\delta
x^h}+XXXY^h_{kij}\frac{\partial}{\partial y^h}
$$
$$
K\big(\frac{\delta}{\delta x^i},\frac{\delta}{\delta x^j}
\big)\frac{\partial}{\partial
y^k}=XXYX^h_{kij}\frac{\delta}{\delta
x^h}+XXYY^h_{kij}\frac{\partial}{\partial y^h}
$$
$$
K\big(\frac{\partial}{\partial y^i},\frac{\partial}{\partial y^j}
\big)\frac{\delta}{\delta x^k}=YYXX^h_{kij}\frac{\delta}{\delta
x^h}+YYXY^h_{kij}\frac{\partial}{\partial y^h}
$$
$$
K\big(\frac{\partial}{\partial y^i},\frac{\partial}{\partial y^j}
\big)\frac{\partial}{\partial
y^k}=YYYX^h_{kij}\frac{\delta}{\delta
x^h}+YYYY^h_{kij}\frac{\partial}{\partial y^h}
$$
$$
K\big(\frac{\partial}{\partial y^i},\frac{\delta}{\delta x^j}
\big)\frac{\delta}{\delta x^k}=YXXX^h_{kij}\frac{\delta}{\delta
x^h}+YXXY^h_{kij}\frac{\partial}{\partial y^h}
$$
$$
K\big(\frac{\partial}{\partial y^i},\frac{\delta}{\delta x^j}
\big)\frac{\partial}{\partial
y^k}=YXYX^h_{kij}\frac{\delta}{\delta
x^h}+YXYY^h_{kij}\frac{\partial}{\partial y^h},
$$
where the $M$-tensor fields appearing as coefficients denote the
horizontal and vertical components of the curvature tensor of the
tangent bundle, and they are given by
$$
XXXX^h_{kij}=\widetilde{S}^h_{il}\widetilde{S}^l_{jk}+P^h_{li}S^l_{jk}
-\widetilde{S}^h_{jl}\widetilde{S}^l_{ik}-P^h_{lj}S^l_{ik}+R^h_{kij}+R_{0ij}^lP^h_{lk}
$$
$$
XXXY^h_{kij}=\widetilde{S}^l_{jk}S^h_{il}+\widetilde{P}^h_{li}S^l_{jk}-
\widetilde{S}^l_{ik}S^h_{jl}-\widetilde{P}^h_{lj}S^l_{ik}+\widetilde{P}^h_{lk}R^l_{0ij}-
$$
$$
-\frac{1}{2}\dot{\nabla}_iR_{0jk}^rG^{(2)}_{rl}H^{(3)}_{hl}+c_3\dot{\nabla}_iR_{j0kh}
$$
$$
XXYX^h_{kij}=\widetilde{P}^l_{kj}P^h_{li}+P^l_{kj}\widetilde{S}^h_{il}-
\widetilde{P}^l_{ki}P^h_{lj}-P^l_{ki}\widetilde{S}^h_{jl}+R^l_{0ij}\widetilde{Q}^h_{lk}
$$
$$
XXYY^h_{kij}=\widetilde{P}^l_{kj}\widetilde{P}^h_{li}+P^l_{kj}S^h_{il}-
\widetilde{P}^l_{ki}\widetilde{P}^h_{lj}-P^l_{ki}S^h_{jl}+R^l_{0ij}Q^h_{lk}+R^h_{kij}
$$
$$
YYXX^h_{kij}=\partial_iP^h_{jk}-\partial_jP^h_{ik}+\widetilde{P}^l_{jk}\widetilde{Q}^h_{il}+
P^l_{jk}P^h_{il}-\widetilde{P}^l_{ik}\widetilde{Q}^h_{jl}-P^l_{ik}P^h_{jl}
$$
$$
YYXY^h_{kij}=\partial_i\widetilde{P}^h_{jk}-\partial_j\widetilde{P}^h_{ik}+\widetilde{P}^l_{jk}Q^h_{il}
+P^l_{jk}\widetilde{P}^h_{il}-\widetilde{P}^l_{ik}Q^h_{jl}-P^l_{ik}\widetilde{P}^h_{jl}
$$
$$
YYYX^h_{kij}=\partial_i\widetilde{Q}^h_{jk}-\partial_j\widetilde{Q}^h_{ik}+
Q^l_{jk}\widetilde{Q}^h_{il}+\widetilde{Q}^l_{jk}P^h_{il}-Q^l_{ik}\widetilde{Q}^h_{jl}
-\widetilde{Q}^l_{ik}P^h_{jl}
$$
$$
YYYY^h_{kij}=\partial_iQ^h_{jk}-\partial_jQ^h_{ik}+Q^l_{jk}Q^h_{il}+\widetilde{Q}^l_{jk}\widetilde{P}^h_{il}-
Q^l_{ik}Q^h_{jl}-\widetilde{Q}^l_{ik}\widetilde{P}^h_{jl}
$$
$$
YXXX^h_{kij}=\partial_i\widetilde{S}^h_{jk}+S^l_{jk}\widetilde{Q}^h_{il}+\widetilde{S}^l_{jk}P^h_{il}
-\widetilde{P}^l_{ik}P^h_{lj}-P^l_{ik}\widetilde{S}^h_{jl}-\dot{\nabla}_jR_{0ik}^rG^{(2)}_{rl}H^{(3)}_{hl}
$$
$$
YXXY^h_{kij}=\partial_iS^h_{jk}+S^l_{jk}Q^h_{il}+\widetilde{S}^l_{jk}\widetilde{P}^h_{il}
-\widetilde{P}^l_{ik}\widetilde{P}^h_{lj}-P^l_{ik}S^h_{jl}-\dot{\nabla}_jR_{0ik}^rG^{(2)}_{rl}H^{(1)}_{hl}
$$
$$
YXYX^h_{kij}=\partial_iP^h_{kj}+\widetilde{P}^l_{kj}\widetilde{Q}^h_{il}+P^l_{kj}P^h_{il}-
Q^l_{ik}P^h_{lj}-\widetilde{Q}^l_{ik}\widetilde{S}^h_{jl}
$$
$$
YXYY^h_{kij}=\partial_i\widetilde{P}^h_{kj}+\widetilde{P}^l_{kj}Q^h_{il}+P^l_{kj}\widetilde{P}^h_{il}-
Q^l_{ik}\widetilde{P}^h_{lj}-\widetilde{Q}^l_{ik}S^h_{jl}.
$$
We mention that we used the character $X$ on a certain position to
indicate that the argument on that position was a horizontal
vector field and, similarly, we used the character $Y$ for
vertical vector fields. \vskip1mm We compute the partial
derivatives with respect to the tangential coordinates $y^i$ of
the entries of the matrices $G$ and $H$
$$
\partial_iG^{(\alpha)}_{jk}=c_\alpha'g_{0i}g_{jk}+d_\alpha'g_{0i}g_{0j}g_{0k}+
d_\alpha g_{ij}g_{0k}+d_\alpha g_{0i}g_{jk}
$$

$$
\partial_iH_{(\alpha)}^{jk}=p_\alpha'g^{jk}g_{0i}+q_\alpha'g_{0i}y^jy^k+
q_\alpha\delta^j_iy^k+q_\alpha y^j\delta^k_i
$$

\vskip1mm
$$
\partial_i\partial_jG^{(\alpha)}_{kl}=c_\alpha''g_{0i}g_{0j}g_{kl}+c_\alpha'g_{ij}g_{kl}+
d_\alpha''g_{0j}g_{0k}g_{0l}+ d_\alpha'g_{ij}g_{0k}g_{0l}+
$$

$$
+d_\alpha'g_{0j}g_{ik}g_{0l}+
d_\alpha'g_{0j}g_{0k}g_{il}+d_\alpha'g_{0i}g_{jk}g_{0l}+d_\alpha'g_{0i}g_{0k}g_{jl}+
d_\alpha g_{jk}g_{il}+d_\alpha g_{ik}g_{jl},
$$
$$
\alpha =1,2,3.
$$
\vskip1mm Next we get the first order partial derivatives with
respect to the tangential coordinates $y^i$ of the $M$-tensor
fields $P^h_{ij},Q^h_{ij},S^h_{ij}, \widetilde P^h_{ij},
\widetilde Q^h_{ij}, \widetilde S^h_{ij}$
$$
\partial_iQ^h_{jk}=\frac{1}{2}\partial_iH_{(2)}^{hl}(\partial_jG^{(2)}_{kl}+\partial_kG^{(2)}_{jl}-
\partial_lG^{(2)}_{jk})+
\frac{1}{2}H_{(2)}^{hl}(\partial_i\partial_jG^{(2)}_{kl}+\partial_i\partial_kG^{(2)}_{jl}
-\partial_i\partial_lG^{(2)}_{jk})+
$$
$$
+\frac{1}{2}\partial_iH_{(3)}^{hl}(\partial_jG^{(3)}_{kl}+\partial_kG^{(3)}_{jl})+
\frac{1}{2}H_{(3)}^{hl}(\partial_i\partial_jG^{(3)}_{kl}+\partial_i\partial_kG^{(3)}_{jl})
$$
$$
\partial_i\widetilde{Q}^h_{jk}=\frac{1}{2}\partial_iH_{(3)}^{hl}(\partial_jG^{(2)}_{kl}+\partial_kG^{(2)}_{jl}-
\partial_lG^{(2)}_{jk})+
\frac{1}{2}H_{(3)}^{hl}(\partial_i\partial_jG^{(2)}_{kl}+\partial_i\partial_kG^{(2)}_{jl}
-\partial_i\partial_lG^{(2)}_{jk})+
$$
$$
+\frac{1}{2}\partial_iH_{(1)}^{hl}(\partial_jG^{(3)}_{kl}+\partial_kG^{(3)}_{jl})+
\frac{1}{2}H_{(1)}^{hl}(\partial_i\partial_jG^{(3)}_{kl}+\partial_i\partial_kG^{(3)}_{jl})
$$
$$
\partial_i\widetilde{P}^h_{jk}=\frac{1}{2}\partial_iH_{(2)}^{hl}(\partial_jG^{(3)}_{kl}-
\partial_lG^{(3)}_{jk})+
\frac{1}{2}H_{(2)}^{hl}(\partial_i\partial_jG^{(3)}_{kl}
-\partial_i\partial_lG^{(3)}_{jk})+
$$
$$
+\frac{1}{2}\partial_iH_{(3)}^{hl}(\partial_jG^{(1)}_{kl}+R^r_{0kl}G^{(2)}_{rj})+
\frac{1}{2}H_{(3)}^{hl}(\partial_i\partial_jG^{(1)}_{kl}+R^r_{ikl}G^{(2)}_{rj}+R_{0kl}^r\partial_iG^{(2)}_{rj})
$$
$$
\partial_iP^h_{jk}=\frac{1}{2}\partial_iH_{(3)}^{hl}(\partial_jG^{(3)}_{kl}-
\partial_lG^{(3)}_{jk})+
\frac{1}{2}H_{(3)}^{hl}(\partial_i\partial_jG^{(3)}_{kl}
-\partial_i\partial_lG^{(3)}_{jk})+
$$
$$
+\frac{1}{2}\partial_iH_{(1)}^{hl}(\partial_jG^{(1)}_{kl}+R^r_{0kl}G^{(2)}_{rj})+
\frac{1}{2}H_{(1)}^{hl}(\partial_i\partial_jG^{(1)}_{kl}+R^r_{ikl}G^{(2)}_{rj}+R_{0kl}^r\partial_iG^{(2)}_{rj})
$$
$$
\partial_iS^h_{jk}=-\frac{1}{2}\{(\partial_i\partial_rG^{(1)}_{jk}+R^l_{ijk}G^{(2)}_{lr}+
R_{0jk}^l\partial_iG^{(2)}_{lr})H_{(2)}^{rh}+
$$
$$
+(\partial_rG^{(1)}_{jk}+R_{0jk}^lG^{(2)}_{lr})\partial_iH_{(2)}^{rh}\}
+c_3'g_{0i}R_{j0kr}H_{(3)}^{rh}+c_3(R_{jikr}H_{(3)}^{rh}+R_{j0kr}\partial_iH_{(3)}^{rh})
$$
$$
\partial_i\widetilde{S}^h_{jk}=-\frac{1}{2}\{(\partial_i\partial_rG^{(1)}_{jk}+R^l_{ijk}G^{(2)}_{lr}
+R_{0jk}^l\partial_iG^{(2)}_{lr})H_{(3)}^{rh}+
$$
$$
+(\partial_rG^{(1)}_{jk}+R_{0jk}^lG^{(2)}_{lr})\partial_iH_{(3)}^{rh}\}
+c_3'g_{0i}R_{j0kr}H_{(1)}^{rh}+c_3(R_{jikr}H_{(1)}^{rh}+R_{j0kr}\partial_iH_{(1)}^{rh}).
$$
It was not convenient to think $c_1,c_2,c_3,d_1,d_2,d_3$ and
$p_1,p_2,p_3$, $q_1,q_2,q_3$ as functions of $t$ since RICCI did
not make some useful factorizations after the command
TensorSimplify. We decided to consider these functions as well as
their derivatives of first, second and  third order, as constants,
the tangent vector $y$ as a first order tensor, the components
$G^{(1)}_{ij}, G^{(2)}_{ij}, G^{(3)}_{ij},$ $H_{(1)}^{ij},
H_{(2)}^{ij}, H_{(3)}^{ij}$ as second order tensors and so on, on
the  Riemannian manifold $M$, the associated indices being
$h,i,j,k,l,r,s.$

\vskip5mm We consider an $n$-dimensional Riemannian manifold $M$
with the fundamental metric $g$. The change of the metric
$$
g^*=\rho^2 g,
$$
where $\rho$ is a certain positive function, does not change the
angle between two vectors at a potint and so is called a
\emph{conformal transformation of the metric.}

The \emph{Weyl conformal curvature tensor} is a tensor field
invariant under any conformal transformation of the metric and it
is given by the expression
$$
C(X,Y)Z=K(X,Y)Z+L(Y,Z)X-L(X,Z)Y+g(Y,Z)NX-g(X,Z)NY,
$$
with $g(NX,Y)=L(X,Y),$ for any vector fields $X,Y,Z$, where $L$ is
a $(2,0)$-tensor field, called by some matematicians \emph{the
tensor of Brinkmann}, given by
$$
L(X,Y)=-\frac{1}{n-2}R(X,Y)+\frac{1}{2(n-1)(n-2)}r g(X,Y).
$$

In local coordinates we have the expressions
$$
C^h_{kij}=K^h_{kji}+\delta^h_kL_{ji}-\delta^h_jL_{ki}+L^h_kg_{ji}-L^h_jg_{ki}
$$
$$
L_{ji}=-\frac{1}{n-2}R_{ij}+\frac{1}{2(n-1)(n-2)}rg_{ji}
$$
$$
L^h_k=L_{kt}g^{th}.
$$

The tensor $C$ vanishes identically for $n=3$.

\vskip2mm If a Riemannian metric $g$ is conformally related to a
Riemannian metric $g^*$ which is locally Euclidian, then the
Riemannian manifold with the metric $g$ is said to be
\emph{conformally flat}. \vskip2mm

For the tangent bundle $TM$ of an $n$-dimensional Riemannian
manifold, with the general natural lifted metric $G$, the
expression of the tensor of conformal curvature becomes
$$
C(X,Y)Z=K(X,Y)Z+L(Y,Z)X-L(X,Z)Y+G(Y,Z)NX-G(X,Z)NY,
$$
where
$$
L(X,Y)=-\frac{1}{2(n-1)}Ric(X,Y)+\frac{1}{4(n-1)(2n-1)}scal\
G(X,Y)
$$
$$
scal=G^{ji}R_{ji};\qquad G(NX,Y)=L(X,Y).
$$

By using the local addapted frame, we obtain
$$
C\big(\frac{\delta}{\delta x^i},\frac{\delta}{\delta x^j}
\big)\frac{\delta}{\delta x^k}=CXXXX^h_{kij}\frac{\delta}{\delta
x^h}+CXXXY^h_{kij}\frac{\partial}{\partial y^h}
$$
$$
C\big(\frac{\delta}{\delta x^i},\frac{\delta}{\delta x^j}
\big)\frac{\partial}{\partial
y^k}=CXXYX^h_{kij}\frac{\delta}{\delta
x^h}+CXYY^h_{kij}\frac{\partial}{\partial y^h}
$$
$$
C\big(\frac{\partial}{\partial y^i},\frac{\partial}{\partial y^j}
\big)\frac{\delta}{\delta x^k}=CYYXX^h_{kij}\frac{\delta}{\delta
x^h}+CYYXY^h_{kij}\frac{\partial}{\partial y^h}
$$
$$
C\big(\frac{\partial}{\partial y^i},\frac{\partial}{\partial y^j}
\big)\frac{\partial}{\partial
y^k}=CYYYX^h_{kij}\frac{\delta}{\delta
x^h}+CYYYY^h_{kij}\frac{\partial}{\partial y^h}
$$
$$
C\big(\frac{\partial}{\partial y^i},\frac{\delta}{\delta x^j}
\big)\frac{\delta}{\delta x^k}=CYXXX^h_{kij}\frac{\delta}{\delta
x^h}+CYXXY^h_{kij}\frac{\partial}{\partial y^h}
$$
$$
C\big(\frac{\partial}{\partial y^i},\frac{\delta}{\delta x^j}
\big)\frac{\partial}{\partial
y^k}=CYXYX^h_{kij}\frac{\delta}{\delta
x^h}+CYXYY^h_{kij}\frac{\partial}{\partial y^h},
$$
where the horizontal and vertical components of the conformal
curvature, are given by the expressions
$$
CXXXX^h_{kij}=XXXX^h_{kij}+LXX_{jk}\delta^h_i-LXX_{ik}\delta^h_j+G^{(1)}_{jk}NXX^h_i-G^{(1)}_{ik}NXY^h_j
$$
$$
CXXXY^h_{kij}=XXXY^h_{kij}+G^{(1)}_{jk}NXY^h_i-G^{(1)}_{ik}NXY^h_j
$$
$$
CXXYX^h_{kij}=XXYX^h_{kij}+LXY_{jk}\delta^h_i-LXY_{ik}\delta^h_j+G^{(3)}_{jk}NXX^h_i-G^{(3)}_{ik}NXX^h_j
$$
$$
CXXYY^h_{kij}=XXYY^h_{kij}+G^{(3)}_{jk}NXY^h_i-G^{(3)}_{ik}NXY^h_j
$$
$$
CYXXX^h_{kij}=YXXX^h_{kij}-LYX_{ik}\delta^h_j+G^{(1)}_{jk}NYX^h_i-G^{(3)}_{ik}NXX^h_j
$$
$$
CYXXY^h_{kij}=YXXY^h_{kij}+LXX_{jk}\delta^h_i+G^{(1)}_{jk}NYY^h_i-G^{(3)}_{ik}NXY^h_j
$$
$$
CYXYX^h_{kij}=YXYX^h_{kij}-LYY_{ik}\delta^h_j+G^{(3)}_{jk}NYX^h_i-G^{(2)}_{ik}NXX^h_j
$$
$$
CYXYY^h_{kij}=YXYY^h_{kij}+LXY_{jk}\delta^h_i+G^{(3)}_{jk}NYY^h_i-G^{(2)}_{ik}NXX^h_j
$$
$$
CYYXX^h_{kij}=YYXX^h_{kij}+G^{(3)}_{jk}NYX^h_i-G^{(3)}_{ik}NYX^h_j
$$
$$
CYYXY^h_{kij}=YYXY^h_{kij}+LYX_{jk}\delta^h_i-LYX_{ik}\delta^h_j+G^{(3)}_{jk}NYY^h_i-G^{(3)}_{ik}NYY^h_j$$
$$
CYYYX^h_{kij}=YYYX^h_{kij}+G^{(2)}_{jk}NYX^h_i-G^{(2)}_{ik}NYX^h_j
$$
$$
CYYYY^h_{kij}=YYYY^h_{kij}+LYY_{jk}\delta^h_i-LYY_{ik}\delta^h_j+G^{(2)}_{jk}NYY^h_i-G^{(2)}_{ik}NYY^h_j,
$$
the horizontal and vertical components of the tensors $L$ and $N$
being
$$
LXX_{ij}=-\frac{1}{2(n-1)}\big(RicXX_{ij}-\frac{1}{2(2n-1)}scal\
G^{(1)}_{ij}\big)
$$
$$
LXY_{ij}=-\frac{1}{2(n-1)}\big(RicXY_{ij}-\frac{1}{2(2n-1)}scal\
G^{(3)}_{ij}\big)
$$
$$
LYX_{ij}=-\frac{1}{2(n-1)}\big(RicYX_{ij}-\frac{1}{2(2n-1)}scal\
G^{(3)}_{ij}\big)
$$
$$
LYY_{ij}=-\frac{1}{2(n-1)}\big(RicYY_{ij}-\frac{1}{2(2n-1)}scal\
G^{(2)}_{ij}\big)
$$
$$
NXX^h_j=LXX_{jk}H_{(1)}^{kh}+LXY_{jk}H_{(3)}^{kh}
$$
$$
NXY^h_j=LXX_{jk}H_{(3)}^{kh}+LXY_{jk}H_{(2)}^{kh}
$$
$$
NYX^h_j=LYX_{jk}H^{kh}_{(1)}+LYY_{jk}H_{(3)}^{kh}
$$
$$
NYY^h_j=LYX_{jk}H^{kh}_{(3)}+LYY_{jk}H_{(2)}^{kh}.
$$

In order to get the conditions under which $(TM,G)$ is a
conformally flat Riemannian manifold, we study the vanishing of
the components of the Weyl conformal curvature tensor. In this
study it is useful the following generic result similar to the
lemma \ref{lema1}.

\begin{lemma}\label{lema2}
If $\alpha _1,\dots , \alpha_{10}$ are smooth functions on $TM$
such that
$$
\alpha_1 \delta^h_i g_{jk}+\alpha_2 \delta^h_j g_{ik}+ \alpha_3
\delta^h_kg_{ij}+\alpha_4 \delta^h_k g_{0i}g_{0j} +\alpha_5
\delta^h_j g_{0i}g_{0k}+\alpha_6 \delta^h_i g_{0j}g_{0k}+
$$
$$
+\alpha_7g_{jk} g_{0i}y^h+ \alpha_8 g_{ik}g_{0j}y^h+\alpha_9
g_{ij}g_{0k}y^h+\alpha_{10}g_{0i}g_{0j}g_{0k}y^h=0,
$$
then $\alpha_1=\dots =\alpha_{10}=0$.
\end{lemma}

After the analysis of the values in $y=0$ of several components of
the Weyl tensor of conformal curvature, computed by using the
RICCI package from Mathematica, we can formulate the next theorem.

\begin{theorem}\label{vanish of R}
Let $(M,g)$ be a Riemannian manifold. If the tangent bundle $TM$
with the general natural lifted metric $G$ is conformally flat,
then the base manifold is of constant sectional curvature.
\end{theorem}
\emph{Proof:} From the vanishing condition for the component
$CXXXX^h_{kij}$ computed in $y=0$, we find the expression of the
Ricci tensor of the base manifold and if we replace this
expression in the component $CXXXY^h_{kij}$ computed in $y=0,$ we
obtain that the Riemannian curvature of the base manifold is given
by
$$
\Rim^h_{kij}= c(g_{jk}\delta^h_i - g_{ik}\delta^h_j) ,
$$
so, the base manifold is of constant sectional curvature.

\section{Conformally flat tangent bundles}

A detailed analysis of the annulation of all the components of the
conformal curvature computed with RICCI, leads to the study of
several cases, which give rise to the next theorems.

\begin{theorem}\label{teorfinala1}
Let $(M,g)$ be a Riemannian manifold and let $G$ be the general
natural lifted metric to the tangent bundle, given by the relations
$(\ref{defG})$. Assume that $c_2+2td_2\neq 0$, $c_3+2td_3\neq 0$,
and
$c_1c_2-c_3^2+2cc_2^2t+2c_1d_2t-4c_3d_3t+4cc_2d_2t^2-4d_3^2t^2\neq
0$. Then the Riemannian manifold $(TM,G)$ is conformally flat if and
only if the base maniflold is flat and the associated matrix of the
natural lifted metric $G$ has one of the forms:
$$
\begin{pmatrix}
 0& \beta g_{ij}+\gamma g_{0i}g_{0j} \\
 \beta g_{ij}+\gamma g_{0i}g_{0j}&\alpha g_{ij}+(\alpha^{'} + \frac{\alpha
 (\gamma-\beta^{'})+ 2\alpha^{'}\gamma t}{\beta} - \frac{2\alpha \beta^{'}\gamma t}{\beta^2})\ g_{0i}g_{0j}
\end{pmatrix}
$$
$$
\begin{pmatrix}
k\ g_{ij}& \beta g_{ij}+\beta^{'}g_{0i}g_{0j} \\
 \beta g_{ij}+\beta^{'}g_{0i}g_{0j}&\alpha g_{ij}+
 \frac{k\alpha^{'}(2 \alpha + {\alpha^{'}}t) -
 2\alpha^{'}\beta(\beta + 2\beta^{'}t) + 4\alpha {\beta^{'}}^2t}{2(k\alpha - \beta^2)}
 \ g_{0i}g_{0j}
\end{pmatrix},
$$
where $k$ is a nonzero arbitrary real constant,
$\alpha,\beta,\gamma$ are some arbitrary real smooth functions
depending on the energy density,  $\alpha\neq \frac{\beta^2}{k}$
and $\beta$ is nonnull.
\end{theorem}

\emph{Proof:} In the proposition \ref{vanish of R} we prooved that
the base manifold of the conformally flat tangent bundle must have
constant sectional curvature, $c$. By using the RICCI package of the
program Mathematica, we replace the corresponding expressions of the
components of $K$ in all the components of the Weyl conformal
curvature tensor of $TM$. After a quite long computation we find
some components in which the third terms are of one of the forms :
$$
\frac{c_1c_3(cc_2-d_1)}{4(c_1c_2-c_3^2))}\delta^h_ky_iy_j\qquad
\rm{in} \qquad CYXXX^{h}_{kij}
$$
$$
\frac{c_3(cc_2-d_1)}{2
(-(c_1c_2+c_3^2))}g_{ij}\delta^h_kg_{ij}\qquad \rm{in} \qquad
CXXXY^{h}_{kij}
$$
$$
\frac{c_1(cc_2-d_1)}{2 (c_1c_2-c_3^2))}g_{ij}\delta^h_k \qquad
\rm{in} \qquad CYXXY^{h}_{kij}
$$
$$
\frac{c_2(cc_2-d_1)}{2 (c_1c_2-c_3^2))}g_{ij}\delta^h_k \qquad
\rm{in} \qquad CYXYX^{h}_{kij}
$$
$$
\frac{c_3(cc_2-d1)}{(2(c_3^2-c_1c_2))}g_{ij}\delta^h_k \qquad
\rm{in} \qquad CYXYY^h_{kij}.
$$

Since the metric $G$ must be non-degenerate,  $c_1$, $c_2$ and
$c_3$ cannot vanish at the same time, so we must have $d_1=cc_2$.
Replacing this expression of $d_1$ in all the components of the
Weyl conformal curvature tensor we get some simpler expressions.

The annulation of the first two coefficients which appear in the
new expression of $CXXXY^h_{kij}$ implies the annulation of the
product
$$
4c(c_1c_2 - c_3^2)^2t(c_1^{'}c_3 + 2cc_2c_3 + 4cc_3d_2t +
2c_1^{'}d_3t),
$$
which leads to two cases: $c=0$ or $c_1^{'}=\frac{2cc_3(c_2 +
2d_2t)}{c_3 + 2d_3t}$.

The above theorem refers to the first case, namely at the case
when the base manifold is flat. The second case will be discussed
in the theorem \ref{teoremafinala3}.

In the first case, from the vanishing condition for the third term
of the component $CXXXX^h_{kij}$, we obtain that
$$
c_1^{'}=\frac{c_1}{c_3}(c_3^{'}-d_3)
$$
and the only coefficient that remains in $CXXXY^h_{kij}$  and
which must vanish is
$$
-c_1^3(c_3^{'}-d_3)(c_3+2td_3)=0.
$$

The annulation of $c_1$ leads, after quite long computations, to
the first form of the matrix and the condition $d_3=c_3^{'}$ leads
to the second form presented in the theorem.

The subcase $c_3+2td_3=0$ will be treated separately (it has been
excluded above), because it is a singular case, and the form
obtained for the associated matrix in this case will be presented in
the theorem \ref{teoremafinala2}. \vskip3mm \textbf{Remark(see
\cite{Druta}):} If the tangent bundle $(TM,G)$ of a Riemannian
manifold $(M,g)$ is of constant sectional curvature, that means the
associated matrix of $G$ has the form
$$
\begin{pmatrix}
k g_{ij} & \beta g_{ij}+
\beta' g_{0i} g_{0j} \\
\beta g_{ij}+ \beta' g_{0i} g_{0j} &\alpha
g_{ij}+\frac{\alpha'\beta^2+2\alpha'\beta
\beta't-2\alpha\beta'^2t}{\beta^2} g_{0i} g_{0j}
\end{pmatrix},
$$
where $\alpha,\beta$ are some arbitrary real smooth functions
depending on the energy density, $\beta$ nonzero, and k is an
arbitrary real constant, then the tangent bundle is conformally
flat if the constant k becomes zero.

\begin{theorem}\label{teoremafinala2}
Let $(M,g)$ be an $n$-dimensional Riemannian manifold and let $G$
be the general natural lifted metric to $TM$, having the
associated matrix of the form obtained in the singular case
$c_3+2td_3=0$, namely
$$
\begin{pmatrix}
 c_1g_{ij}+
d_1g_{0i}g_{0j}& -2td_3g_{ij}+
d_3g_{0i}g_{0j} \\
 -2td_3g_{ij}+
d_3g_{0i}g_{0j}&c_2g_{ij}+ d_2g_{0i}g_{0j}
\end{pmatrix},
$$
where $c_1,d_1,c_2,d_2,d_3$ are smooth real functions depending on
the density of energy. Assume that $c_1^4c_2^2 + 2c_1^4c_2c_2't +
c_1^2c_1'^2c_2^2t^2 + c_1^4c_2'^2t^2 - 32c_1'^2d_3^4t^6\neq 0$.
Then the bundle of nonzero tangent vectors to $M$, $TM_0$, is
conformally flat with respect to the natural lifted metric $G$, if
and only if the base manifold is flat and the matrix becomes
$$
\begin{pmatrix}
 kg_{ij}& \frac{e^{\varepsilon}}{t\sqrt{t}}(-2tg_{ij}+g_{0i}g_{0j}) \\
\frac{e^{\varepsilon}}{t\sqrt{t}}(-2tg_{ij}+g_{0i}g_{0j})&\alpha
g_{ij}+ \frac{k\alpha^{'}(2\alpha + \alpha^{'}t) +
4\alpha\frac{e^{2{\varepsilon}}}{t^2}}{2(k\alpha -
\frac{4e^{{2\varepsilon}}}{t})}g_{0i}g_{0j}
\end{pmatrix},
$$
where $k$ and $\varepsilon$ are two arbitrary real constants,
$k\neq 0$ and $\alpha$ is a smooth real function depending on the
density of energy, $\alpha\neq \frac{4e^{2\varepsilon}}{kt}$.
\end{theorem}

\begin{corollary}
Let $(M,g)$ be an $n$-dimensional Riemannian manifold. The bundle
of nonzero tangent vectors to $M$, $TM_0$, is conformally flat
with respect to the natural lifted metric $G$ of diagonal type, if
and only if the base manifold is flat, and the associated matrix
of $G$ has the form
$$
\begin{pmatrix}
kg_{ij}
& 0 \\
 0&\alpha g_{ij}+ \frac{k\alpha^{'}(2\alpha + \alpha^{'}t) +
4\alpha\frac{e^{2\varepsilon}}{t^2}}{2(k\alpha -
\frac{4e^{2\varepsilon}}{t})}g_{0i}g_{0j}
\end{pmatrix},
$$
where $k$ and $\varepsilon$ are two arbitrary real constants,
$k\neq 0$ and $\alpha$ is a smooth real function depending on the
density of energy, $\alpha\neq \frac{4e^{2\varepsilon}}{kt}$.
\end{corollary}

The next theorem presents the form of the matrix associated to the
general natural lifted metric $G$, obtained in the second case
mentioned in the proof of the theorem \ref{teorfinala1}, namely
$c_1^{'}=\frac{2cc_3(c_2 + 2d_2t)}{c_3 + 2d_3t}$. By using RICCI,
we obtain that all its subcases reduce to some cases that we have
already treated, except the singular subcase $c_2+2td_2=0$. So,
the next theorem can be formulated as follows:

\begin{theorem}\label{teoremafinala3}
Let $(M,g)$ be a Riemannian manifold, and let $G$ be the natural
lifted metric to $TM$, having the associated matrix of the form
obtained in the singular case $c_2+2td_2=0,$ namely
$$
\begin{pmatrix}
c_1g_{ij}+d_1g_{0i}g_{0j}& c_3g_{ij}+d_3g_{0i}g_{0j}\\
c_3g_{ij}+d_3g_{0i}g_{0j}&-2td_2g_{ij}+d_2g_{0i}g_{0j}
\end{pmatrix},
$$
where $c_1, d_1, c_2, d_2, c_3,d _3$ are some arbitrary smooth
real functions of the energy density. The tangent bundle $TM$ is
conformally flat with respect to the natural metric $G$ if and
only if the base manifold is flat and the matrix associated to $G$
becomes of the antidiagonal form
$$
\begin{pmatrix}
0& c_3g_{ij}+d_3g_{0i}g_{0j}\\
c_3g_{ij}+d_3g_{0i}g_{0j}&0
\end{pmatrix}.
$$
\end{theorem}

\vskip3mm

S.L. Dru\c t\u a

Faculty of Mathematics

"Al.I. Cuza" University of Ia\c si

Bd. Carol I, no. 11

700506 Ia\c si

ROMANIA

\email{simonadruta@yahoo.com}


\begin{thebibliography}{99}

\bibitem{Abbassi1} {\bf  Abbassi, M.T.K.,  Sarih, M.,}  {\it On some hereditary propertoes of
Riemannian g-natural metrics on tangent bundles of Riemannian
manifolds,} Diff. Geom. And its Appl. 22 (2005), 19-47.

\bibitem{Abbassi2}{\bf Abbassi, M.T.K., Sarih, M.,} {\it On Riemannian g-natural metrics of the
form $a.g^s+b.g^h+c.g^v $ on the tangent bundle of a Riemannian
manifold (M,g)}. Mediterranean J. Math. 2 (2005), 19-43.
\bibitem{BejanOpr} {\bf  Bejan, C. L., Oproiu, V.,} {\it Tangent bundles of
quasi-constant holomorphic sectional curvatures},  Balkan J. Geom.
Applic., 11 (2006), 11-22 ]
\bibitem{Druta} {\bf Druta, S. L.,} {\it The sectional curvature
of tangent bundles with general natural lifted metrics}, to appear
in Proceedings of the Ninth International Conference on Geometry,
Integrability and Quantisation, 8-13 June, 2007,  Varna, Bulgaria,
Eds. I. M. Mladenov and M. Leon, Sofia.
\bibitem{Kolar} {\bf Kol\'a\v r,~I., Michor,~P., Slovak,~J.,} {\it Natural
Operations in Differential Geometry}, Springer Verlag, Berlin,
1993, vi, 434 pp.


\bibitem{KowalskiSek} {\bf Kowalski,~O., Sekizawa,~M.,} {\it Natural
transformations of Riemannian metrics on manifolds to metrics on
tangent bundles - a classification,} Bull. Tokyo Gakugei Univ.
(4), 40 (1988), 1-29.
\bibitem{Krupka} {\bf Krupka,~D., Jany\v ska,~J.,} {\it Lectures on
Differential Invariants,} Folia Fac. Sci. Nat. Univ. Purkinianae
Brunensis, 1990.
\bibitem{Mok} {\bf Mok,~K.P., Patterson,~E.M., Wong,~Y.C.,} {\it Structure
of symmetric tensors of type (0,2) and tensors of type (1,1) on
the tangent bundle},~ Trans. Amer. Math. Soc.,~234
(1977),~253-278.
\bibitem{Munteanu1} {\bf Munteanu,~M.} {\it Cheeger Gromoll type metrics on the tangent
bundle,} Proceedings of Fifth International Symposium
BioMathsPhys, Iasi , June 16-17, 2006, 9pp.
\bibitem{Munteanu2} {\bf Munteanu,~M.} {\it Old and New Structures on the Tangent
Bundle,} Proceedings of the Eighth International Conference on
Geometry, Integrability and Quantization, June 9-14, 2006, Varna,
Bulgaria, Ed. I. M. Mladenov and M. de Leon, Sofia 2007 264-278.
\bibitem{Oproiu4} {\bf Oproiu,~V.,} {\it A generalization of natural
almost Hermitian structures on the tangent bundles.} Math. J.
Toyama Univ., 22 (1999), 1-14.
\bibitem{Oproiu1} {\bf Oproiu,~V.,} {\it A locally symmetric Kaehler Einstein
structure on the tangent bundle of a space form}, Beitr\"age
Algebra Geom/Contributions to Algebra and Geometry, 40 (1999),
363-372.
\bibitem{Oproiu2} {\bf Oproiu,~V.,} {\it A Kaehler Einstein structure on the
tangent bundle of a space form}, Int. J. Math. Math. Sci. 25 (3)
(2001), 183-195.
\bibitem{Oproiu3} {\bf Oproiu,~V.,} {Some new geometric structures on the
tangent bundles}, Publ. Math. Debrecen, 55/3-4 (1999), 261-281.
\bibitem{OprDruta}{\bf  Oproiu,~V., Druta,~S.,} {\it General natural K\"ahler structures of constant holomorphic
sectional curvature on tangent bundles},~An.St.Univ. "Al.I.Cuza"
Iasi, Matematica,
\bibitem{OprPap1}  {\bf Oproiu,~V., Papaghiuc,~N.,} {\it A Kaehler structure
on the nonzero tangent bundle of a space form,} Diff. Geom. Appl.
11 (1999), 1-14.
\bibitem{OprPap2} {\bf Oproiu,~V., Papaghiuc,~N.} {\it Some classes
of almost anti-Hermitian structures on the tangent bundle},
Mediterranean Journal of Mathematics  1 (3) (2004), 269-282.
\bibitem{Tahara2} {\bf Tahara,~M., Vanhecke,~L., Watanabe,~Y.,} {\it New
structures on tangent bundles}, Note di Matematica (Lecce), 18
(1998), 131-141.
\bibitem{YanoIsh} {\bf Yano,~K., Ishihara,~S.,} {\it Tangent and Cotangent
Bundles},~M. Dekker Inc., New York,~1973.


\end{thebibliography}
\end{document}